\documentclass[12pt]{article}
\newcommand{\bc}{\begin{center}}
\newcommand{\ec}{\end{center}}
\newcommand{\be}{\begin{equation}}
\newcommand{\ee}{\end{equation}}
\newcommand{\bea}{\begin{eqnarray}}
\newcommand{\eea}{\end{eqnarray}}
\newcommand{\ba}{\begin{array}}
\newcommand{\ea}{\end{array}}
\newcommand{\edc}{\end{document}}

\begin{document}
\thispagestyle{empty}
\begin{center}

{\Large {\bf On periodic $p$-harmonic functions on Cayley tree.}}\\
\vspace{0.4cm}
{\bf U.A. Rozikov$^1$, F.T. Ishankulov$^2$ }\\
$^1$ Institute of Mathematics and Informational Technologies,\\
Tashkent,
 Uzbekistan. E-mail: rozikovu@yandex.ru\\
$^2$Samarkand State University, Samarkand, Uzbekistan. E-mail: fandor83@mail.ru\\
 [2mm]
\end{center}
\vspace{0.5cm}

{\bf Abstract:} We show that any periodic with respect to normal
subgroups (of the group representation of the Cayley tree) of
finite index $p$-harmonic function is a constant. For some normal
subgroups of infinite index we describe a class of (non-constant)
periodic $p$-harmonic functions. If $p\neq2$, the $p$-harmonicity
is non-linear, i.e., the linear combination of $p$-harmonic
functions need not be $p$-harmonic. In spite of this, we show that
linear combinations of the $p$-harmonic functions described for
normal subgroups of infinite index are also $p$-harmonic.

\section{Introduction} A tree is connected acyclic graph. One
special case of this is a Cayley tree, i.e., an infinite tree in
which each vertex has exactly $k+1$ incident edges. The Cayley
tree can be represent as the group $G_k$ which is the free product
of $k+1$ second order cyclic groups [1],[2],[6]. The group
representation of the Cayley tree was used in [1],[2],[5]-[8] to
study models of statistical mechanics and describe the sets of
periodic Gibbs measures, also to study random walk trajectories in
a random medium on the Cayley tree. These problems are related to
description of harmonic functions on trees. In [5] a natural
generalization of the notion of harmonic functions on a Cayley
tree is introduced. Using some properties of $G_k$ a set of
harmonic functions described.
    Harmonic functions play an important role in probability
    theory, dynamical systems, statistical mechanics, and theory
    of electrical circuits.
    Note that $p$-harmonic functions are zeros of the
    $p$-Laplacian defined by
$$ \Delta_pu={\rm div}\left(|\nabla u|^{p-2} \nabla u\right), \ \ 1<p<\infty.$$
If $p\neq2$, then $p$-harmonicity is non-linear. We consider the
equation $\Delta_pu=0$ on (discrete) Cayley tree. Since there is
no theory of functional equations with the unknown functions
defined  on a tree, the problem of finding all solution of
$\Delta_pu=0$ on tree is by no means easy. On the other hand, the
theory of functional equations is not very developed even for case
in which the unknown function is defined on $R$, only the
solutions of functional equations of some special forms are known
(e.g. see [3], [4]).
    Thus, it is natural to find periodic (simple) solutions of
    $\Delta_pu=0$ first.

    The main goal of this paper is to describe the set of
    $p$-harmonic functions which are periodic with respect to
    subgroups of the group $G_k$.

    In section 2 we give some preliminary definitions and results
    which needed when describing $p$-harmonic (periodic)
    functions. Section 3 devoted to description of periodic
    $p$-harmonic functions (with respect to normal subgroups of
    finite index). We prove that only constant functions can be
    such periodic. Section 4 devoted to the problem in the case of infinite
    index. We find (explicitly) some $p$-harmonic functions which
    are periodic with respect to a subgroup of infinite index. In
    the last section we show that linear combinations of these
    functions are also $p$-harmonic.
    \section{Definitions and statement of the problem}
    \subsection{Cayley tree.}
    A Cayley tree is an infinite tree in which each vertex has
    exactly $k+1$ incident edges (the Cayley tree of order
    $ k\geq 1$). Let $\Gamma^k=\left(V,L\right)$ be the Cayley
    tree of order $ k\geq 1$, where $V$ and $L$ are the vertex set
    and the edge set, respectively, of $\Gamma^k$.

    If  $x,y\in V$ are the endpoints of an edge $l\in L$, then $x$
    and $y$ are said to be adjacent; in this case, we write $l=<x,y>$.
    The Cayley tree is equipped with a distance $d(x,y),x,y\in V$,
    given by the formula
    \bea
    d(x,y)=\min\{ d: \ \exists x=x_0,x_1,...,x_{d-1},x_d=y\in V
\ \mbox{where} \nonumber \\
 \ <x_0,x_1>,...,<x_{d-1},x_d> \ \mbox{are adjacent} \}.\nonumber
\eea
    The sequence $\pi=\{x=x_0,x_1,...,x_{d-1},x_d=y\in V\}$ at which
    the minimum is attained is called the path form $x$ to $y$.
    Let $G_k$ be the group representation of the Cayley tree;
    i.e., $G_k$ is the free product of $k+1$ second - order
     cyclic groups with generators $a_1,...,a_{k+1}$ such that
     ${{a}}^{2}_i=e,  i=1,...,k+1$,where $e\in G_k$ is the
     identity element.
     Let $x\in G^k$. We denote by
     $$S(x)=\{y\in G_k:<x,y>\}$$
     the set of vertices adjacent to $x$.
\subsection{Normal subgrups of $G_k$.}
    Let $\omega_x(a_i)$ be the number of occurrences of the letter
    $a_i,  i=1,...,k+1$, in the redused word $x\in G_k$.
    Let $A\subseteq N_k=\{1,...,k+1\}$. The set $H_A=\{x\in
    G_k:\sum\limits_{i\in A}\omega_x(a_i)\ \rm is\ \rm even\}$ is a normal
    subgroup of index 2 of $G_k$. Some set of normal subgroups with index $2^m$ can be
    obtained by intersection
    $H_{A_1}\bigcap H_{A_2}\bigcap...\bigcap
    H_{A_m}$ for suitable $A_1,...,A_m\subseteq N_k$ (see. [2]).
    The following theorem gives possible values of the (finite)
    index of normal subgroups.
\vskip 0.5 truecm THEOREM 1. [5]{\it The group $G_k$ does not have
normal subgroups of odd index $(\neq1)$. It has a normal subgroup
of arbitrary even index.}

    Also there are normal subgroups of infinite index. Some of them
    can be described as following [7].
    Fix $M\subseteq N_k$ such that $|M|>1$. $|\bullet|$ is the
    cordinality of $\bullet$.
    Let the mapping $\pi_M:\{a_1,...,a_{k+1}\}\longrightarrow \{a_i,
    \ i\in M\}\cup \{e\}$ be defined by
    $$\pi_M(a_i)=\left\{%
\begin{array}{ll}
    a_i, & \hbox{if} \ \ i\in M \\
    e, & \hbox{if} \ \ i\notin M. \\
\end{array}
\right.$$ Denote by $G_M$ the free product of cyclic groups
$\{e,a_i\}, \ i\in M$. Consider
$$f_M(x)=f_M(a_{i_1}a_{i_2}...a_{i_m})=\pi_M(a_{i_1})
\pi_M(a_{i_2})...\pi_M(a_{i_m}).$$
    Then it is easy to see that $f_M$ is a homomorphism and hence
    $H_M=\{x\in G_k: \ f_M(x)=e\}$ is a normal subgroup of
    infinite index.

\subsection{$p$-Harmonic functions.}
The resistances $r$ is a positive function on $L$ such that
$r(y,x)=r(x,y)$. Denote ${\cal R}=\{r(x,y):(x,y)\in L\}$. We
define the discrete $\nabla u$ and the discrete $p$-Laplacian
$\Delta_p u$ for a function $u$ on $V\left(=G_k\right)$ by$$\nabla
u(x,y)=r(x,y)^{-1}\left(u(y)-u(x)\right),$$ $$\Delta_p
u(x)=\sum\limits_{y\in S(x)}|\nabla u(x,y)|^{p-2}\nabla u(x,y),$$
where $1<p<\infty$. Let $D\subset V$. If $\Delta_p u=0$ in $D$,
then we say that $u$ is $p$-harmonic in $D$ [9]-[11].
    Let $\{u_1,...,u_n\}$ be an $m$-tuple of $p$-harmonic
    functions in $D$. If $p\neq 2$, then the linear combination of
    $p$-harmonic functions $u_1,...,u_m$ need not be $p$-harmonic.
    The $m$-tuple of $p$-harmonic functions in $D$ has a linear
    relation if every linear combination $\sum_{j=1}^m{t_ju_j}$ is
    $p$-harmonic in $D$. Also we say that $\{u_1,...,u_m\}$ has a
    partial linear relation in $D$ if $\sum_{j=1}^m{t_ju_j}$
    $p$-harmonic in $D$ for some $t_1,...,t_m\in R\setminus
    \{0\}$.
    \vskip 0.5 truecm  DEFINITION 1. Let $H \subset G_k$ be a subgroup. A
    function $\{u(x), x\in G_k\}$ is said to be $H$-periodic if
    $u(x)=u(yx)$ for any $x\in G_k$ and $y\in H$.

    If $H$ is a normal subgroup of finite index, then the
    description of $H$-periodic solutions of $\Delta_p u$ can be
    reduced to solving a system of equations with finitely many
    unknowns.

    The main goal of this paper is to describe $H$-periodic $p$-harmonic
    functions for any normal subgroup $H$ and $\forall
    p\in(1,\infty)$.
    \section{The case of finite index.}
     Let $G^{*}_k$ be a
normal subgroup of (finite) index $m\geq 1$ for $G_k$.
    Denote by ${\cal H}_{p,k,m}\left(G^{*}_k, {\cal R}\right)$ the set of all
    $G^{*}_k$-periodic, $p$-harmonic functions for given
    $p\in(1,\infty)$, ${\cal R}$, $m\geq 1$ and $k\geq 1$.
    The following theorem completely describes the set of all
    periodic $p$-harmonic functions for finite $m\geq 1$.
\vskip 0.5 truecm THEOREM 2. {\it For $\forall p\in(1,\infty)$,
$\forall m\geq 1$, $\forall k\geq 1$, $\forall G^{*}_k\subset G_k$
with index $m$ and any fixed ${\cal R}=\{r(x,y)>0:(x,y)\in L\}$
the set ${\cal H}_{p,k,m}\left(G^{*}_k,{\cal R}\right)$ contains
only constant functions.}
 \vskip 0.5 truecm {\it Proof.} Let
 $G_k|G^{*}_k=\{G^{1}_k,...,G^{m}_k\}$ be the factor group. Any
 $G^{*}_k$-periodic $p$-harmonic function $u(x)$ has the
 following form $$u(x)=u_i \ {\rm  if}\ x\in G^{i}_k, i=1,...,m\eqno(1)$$
 where $u_i,\ i=1,...,m$ are a solution to the following system
 $$\sum_{j=1}^m\sum_{y\in S_j(x)}{|u_j-u_i|^{p-2}\over r^{p-1}(x,y)}(u_j-u_i) \mathbf{1}(S_j(x)\neq \emptyset)=0,
 \ \ i=1,...,m\eqno(2)$$ where $x\in G^{i}_k,\ \ S_j(x)=S(x)\cap G^{j}_k, \ \
 j=1,...,m$.
    Note that some $S_j(x)$ can be empty set, so

$$\mathbf{1}(S_j(x)\neq \emptyset)=
\left\{%
\begin{array}{ll}
    1, & \hbox{if} \ \ S_j(x) \neq \emptyset \\
    0, & \hbox{if} \ \ S_j(x)=\emptyset \ . \\
\end{array}
\right.$$
    We shall prove that the system (2) has only solutions
    $(u^{*}_1,...,u^{*}_m)$ with $u^{*}_1=...=u^{*}_m$.
    Assume that there is a solution $u^*=(u^{*}_1,...,u^{*}_m)$
    with $$u^{*}_i\neq u^{*}_j \ \  {\rm for } \  {\rm some} \ \ i,j\in
    \{1,...,m\}\eqno(3)$$
    Denote $u^{*}_{i_0}={\rm max}\{u^{*}_i, \ \ i=1,...,m\}$.

From $i_0$th equation of (2) we have
$$\sum_{j=1}^m\sum_{y\in S_j(x)}{|u^{*}_j-u^{*}_{i_0}|^{p-2}
\over r^{p-1}(x,y)}\mathbf{1}(S_j(x)\neq \emptyset)(u^{*}_j-u^{*}_{i_0})=$$
$$\sum_{j=1:\atop S_j(x)\neq \emptyset}^m\sum_{y\in S_j(x)}{|u^{*}_j-u^{*}_{i_0}|^{p-2}(u^{*}_j-u^{*}_{i_0})\over
r^{p-1}(x,y)}<0,\eqno(4)$$ since $u^{*}_j-u^{*}_{i_0}\leq0$ and by
(3) ${u^{*}_j-u^{*}_{i_0}}<0$ for at least one
$j\in\{1,...,m\}\setminus\{i_0\}$.
    The inequality (4) contradicts to the assumption that $u^{*}$
    is a solution of (2). Theorem is proved.
    \section{The Case of infinite index.}
    In this section we show that for some subgroups of infinite index
there are non - constant $p$-harmonic functions. Consider
$M=\{i,j\}, \ i\neq j\in N_k$ and $H_M=H_{ij}=\{x\in G_k: \
f_M(x)=e\}$ (see section 2.2).
    It is known that $H_{ij}$ is a normal subgroup of infinite
    index and the corresponding factor group can be written as
    (see [7]): $$
    G_k|H_{ij}=\{...,H_{-2},H_{-1},H_0,H_1,H_2,...\}.$$
\vskip 0.5 truecm LEMMA 1. {\it If $x\in H_n$ then $S(x)\subset
H_{n-1}\cup H_n\cup H_{n+1}$. Moreover $|S(x)\cap H_{n-1}|=1, \
|S(x)\cap H_n|=k-1, \ |S(x)\cap H_{n+1}|=1.$}
 \vskip 0.5 truecm {\it Proof.} If $x\in H_n$ then
 $f_M(x)=\underbrace{{a_ia_j \dots a_i}}_{n}$ (or $a_ia_j...a_j$ depending on $n$).
 $S(x)=\{xa_s, \ s=1,...,k+1\}$. We have
 $$f_M(xa_s)=f_M(x)f_M(a_s)=a_ia_j \dots a_if_M(a_s)=\left\{%
\begin{array}{ll}
    \underbrace{{a_ia_j \dots a_j}}_{n-1} & \hbox{if} \ \  a_s=a_i \\
   \underbrace{a_ia_j \dots }_{n} & \hbox{if} \ \ a_s \neq a_i,a_j \\
    \underbrace{a_ia_j \dots a_ia_j}_{n+1}  & \hbox{if} \ \ a_s=a_j \\
\end{array}%
\right.    $$
    This relation complete the proof.

    Note that any $H_{ij}$-periodic function $u(x)$ has the form
    $$u(x)=u_n \ {\rm if} \ x\in H_n \eqno(5)$$
    Consider $H_{ij}$-periodic collection ${\cal R}$ of resistance
    functions $r(x,y)$ i.e. $r(x,y)=r_{nm} \ {\rm if} \  x\in H_n,y\in H_m$.
    By lemma 1 the equation $\Delta_p u=0$ can be written as
    $${{{|u_{n+1}-u_n|^{p-2}\over r^{p-1}_{n,n+1}}(u_{n+1}-u_n)}+{|u_{n-1}-u_n|^{p-2}\over
    r^{p-1}_{n-1,n}}(u_{n-1}-u_n)}=0\eqno(6)$$
    for any $n\in Z$.

    Denote $a_n=u_{n+1}-u_n$ then from (6) we see that
    $a_n,a_{n-1}$ must have the same sign i.e., $a_n\cdot
    a_{n-1}>0, \ n\in Z$.

    Thus $u_n, \ n\in Z$ must be a monotone sequence.

    Put
    $$X_n={|a_n|^{p-2}\over r^{p-1}_{n,n+1}}a_n.$$
    From (6) we get
    $$X_n-X_{n-1}=0, \ \forall n\in Z.\eqno(7)$$
    Hence $X_n=C={\rm Const} \  {\rm for} \ {\rm any} \ n\in Z$.

    Since $u_n$ must be monotone, we have two possibility:

    1) $u_{n+1}>u_n$ then $C>0$. From $X_n=C$ we get
    $$u_{n+1}=C^{1\over{p-1}}\cdot r_{n,n+1}+u_n, \ n\in Z.\eqno(8)$$
    From (8) we obtain
    $$u_n=C^{1\over{p-1}}\sum_{s=-\infty}^{n-1}r_{s,s+1}, \ n\in
    Z.\eqno(9)$$
    2) $u_{n+1}<u_n$ then $-C>0$. In this case we get
 $$u_n=(-C)^{1\over{p-1}}\sum_{s=n}^{+\infty}r_{s,s+1}, \ n\in
    Z.\eqno(10)$$
    Thus we have proved the following
\vskip 0.5 truecm THEOREM 3. {\it Let ${\cal R}$ be
$H_{ij}$-periodic, and
$\sum_{s=-\infty}^{\infty}r_{s,s+1}<+\infty$.
    Then there are two family $U_1,U_2$ of $H_{ij}$-periodic
    $p$-harmonic functions on the Cayley tree of $k\geq1$ such
    that}
    $$U_1=\{u: \
    u(x)=u_n=C^{1\over{p-1}}\sum_{s=-\infty}^{n-1}r_{s,s+1} \ {\rm
    if} \ x\in H_n, \ n\in Z, \ C\geq 0 \};$$
    $$U_2=\{u: \
    u(x)=u_n=C^{1\over{p-1}}\sum_{s=n}^{+\infty}r_{s,s+1} \ {\rm
    if} \ x\in H_n, \ n\in Z, \ C\geq 0 \}.$$
    \vskip 0.5 truecm{\bf Remarks.}

    1. Theorem 3 is also true for
    more general class of ${\cal R}$ than the class of $H_{ij}$-periodic
    resistance functions. For example, one can prove theorem 3 for resistance
    functions:
$$r(x,y)=\left\{%
\begin{array}{ll}
    r_{nm}, & \hbox{if} \ \ x\in H_n, y\in H_m, \ n\neq m \\
    r(x,y), & \hbox{if} \ \  x,y\in H_n \\
\end{array}%
\right.$$ i.e. $r$ takes arbitrary values on $(x,y)\in L$ such
that $x,y\in H_n$.

    2. The role of $p$ in $U_1$ and $U_2$ is not important: one
    can  denote $C^{1\over{p-1}}$ by $C$.

    3. Note that for $M\subset N_k$ the problem of describing of $
H_M$-periodic $p$-harmonic functions on Cayley tree of order $k$
is equivalent to describing of arbitrary (non - periodic)
$p$-harmonic functions on a Cayley tree of order $|M|-1$. Thus for
$M\subset N_k$ with $|M|\geq 3$ this problem is difficult.

    4. By definition a constant is a $p$-harmonic function and the
    linear combination of an arbitrary $p$-harmonic function and a
    constant is also $p$-harmonic. Thus the classes $U_1$ and
    $U_2$ generate two parametrized classes of $p$-harmonic
    functions: $U_1+C_1$, $U_2+C_2$ where $U_i+C_i$ means to each
    functions of $U_i$ a constant $C_i$ is added.
    \section{Linear relations.}
    In this section we describe $p$-harmonic linear combinations of
$p$-harmonic functions from $U_1, \ U_2$. \vskip 0.5 truecm
THEOREM 4. {\it Let $\{v_1,...,v_{q_1},v_{{q_1}+1},...,v_q\}$ be a
$q$-tuple of $p$-harmonic functions with
$\{v_1,...,v_{q_1}\}\subset U_1$ and
$\{v_{{q_1}+1},...,v_q\}\subset U_2$. Then $\{v_1,...,v_q \}$ has
a linear relation on Cayley tree.} \vskip 0.5 truecm {\it Proof.}
Note that $v=\sum_{i=1}^{q}t_iv_j$, $t_i\in R\setminus\{0\}$ has
the following form
$$v(x)=\sum_{i=1}^{q_1}t_i\left(C_i\sum_{s=-\infty}^{n-1}r_{s,s+1}\right)
+\sum_{i=q_1+1}^{q}t_i\left(C_i\sum_{s=n}^{+\infty}r_{s,s+1}\right)\eqno(11)$$
if $x\in H_n, \ n\in Z$. Hence $v$ is a $H_{ij}$-periodic
function. We shall prove that $\Delta_p v=0$ i.e. $v$ is
$p$-harmonic.

    Denote $\varphi_p(t)=|t|^{p-2}t$.
    By definition and (11) we have
    $$\Delta_pv=\sum\limits_{y\in S(x)}{|v(y)-v(x)|^{p-2}\over r^{p-1}(x,y)}(v(y)-v(x))=
    \sum\limits_{y\in S(x)}{\varphi_p(v(y)-v(x))\over
    r^{p-1}(x,y)}=$$
     $${\varphi_p\left(\sum_{i=1}^{q_1}t_iC_i\left(\sum_{s=-\infty}^{n-2}r_{s,s+1}-\sum_{s=-\infty}^{n-1}r_{s,s+1}\right)+
     \sum_{i=q_1+1}^{q}t_iC_i\left(\sum_{s=n-1}^{\infty}r_{s,s+1}-\sum_{s=n}^{\infty}r_{s,s+1}\right)\right)\over r^{p-1}_{n-1,n}}+$$
     $${\varphi_p\left(\sum_{i=1}^{q_1}t_iC_i\left(\sum_{s=-\infty}^{n}r_{s,s+1}-\sum_{s=-\infty}^{n-1}r_{s,s+1}\right)+
     \sum_{i=q_1+1}^{q}t_iC_i\left(\sum_{s=n+1}^{\infty}r_{s,s+1}-\sum_{s=n}^{\infty}r_{s,s+1}\right)\right)\over
     r^{p-1}_{n,n+1}}=$$
     $${{\varphi_p\left(r_{n-1,n}\left(-\sum_{i=1}^{q_1}t_iC_i + \sum_{i=q_1+1}^{q}t_iC_i\right)\right)\over
r^{p-1}_{n-1,n}}+{\varphi_p\left(r_{n+1,n}\left(\sum_{i=1}^{q_1}t_iC_i
- \sum_{i=q_1+1}^{q}t_iC_i\right)\right)\over
r^{p-1}_{n,n+1}}}\eqno(12)$$

    It is easy to see that
    $\varphi_p(st)=\varphi_p(s)\varphi_p(t)$;
    $\varphi_p(t)=t^{p-1} \ {\rm if} \ t>0$ and
    $\varphi_p(-t)=-\varphi_p(t)$. Using these properties from
    (12) we get $\Delta_pv=0$. The theorem is proved.
\vskip 0.5 truecm {\bf Remark.}

The $H_{ij}$-periodic functions described in Theorem 3 and 4 do
not depend on $i,j\in N_k, \ i\neq j$.

 \vskip 0.5 truecm {\bf References}

1. N.N. Ganikhadjaev, The group representation and authomorphisms
of the Cayley tree. Dokl. Akad. Nauk. Resp. Uzb. (1994), ¹4, 3-5.

2. N.N. Ganikhadjaev, U.A. Rozikov. Description of nonperiodic
extreme Gibbs measures of some models on Cayley tree. Theor. Math.
Phys. 111 (1997), ¹1, 109-117.

3.M. Kuczma. Functional equations and their applications, Academic
Press, N.York-London, 1966.

4. M. Kuzma, B. Choczewski, and R.Ger. Iterative functional
equations. Cambridge University Press, Cambridge, 1990.

5. E.P. Normatov, U.A. Rozikov. A description of harmonic
functions via properties of the group representation of the Cayley
tree. Math. Notes. 79 (2006) ¹3, 399-407.

6. U.A. Rozikov. Structure of partitions of the Cayley tree and
their applications to the description of periodic Gibbs
distributions. Theor. Math. Phys. 112 (1997), ¹1, 170-175.

7. U.A. Rozikov. Countably periodic Gibbs measures of the Ising
model on the Cayley tree. Theor. Math. Phys. 130 (2002), ¹1,
92-100.

8. U.A. Rozikov. Random walks in random media on a Cayley tree.
Ukrain. Math. J. 53 (2001), ¹10, 1391-1401.

9. P.M. Soardi. Potential theory of infinite networks, V. 1590.
Lecture Notes in Mathematics. Springer-Verlag, Berlin, 1994.

10. M. Yamasaki. Parabolic and hyperbolic infinite
networks,Hiroshima Math. J. 7 (1977) ¹1, 135-146.

11. M. Yamasaki.
Nonlinear Poisson equations on an infinite network. Mem. Fac. Sci.
Shimane Univ. 23 (1989), 1-9.

\end{document}